\documentclass[12pt]{amsart}

\usepackage{amsmath, amsthm, amscd, amsfonts, amssymb, graphicx, color}
\usepackage{latexsym}
\usepackage{amscd}
\usepackage{amsthm}
\usepackage{amssymb}
\usepackage{enumerate}
\usepackage{mathabx}
\usepackage{color}
\usepackage{soul}
\setstcolor{red}
\setul{}{0.4ex}

\theoremstyle{plain}
\newtheorem{theorem}{Theorem}

\theoremstyle{remark}

\begin{document}

\title{A new proof of a classical result on the topology of orientable connected and compact surfaces by means of the Bochner technique}

\author{J. M.~Almira$^*$, A. Romero}

\thanks{$^*$ Corresponding author. The second author is partially supported by Spanish MINECO and ERDF project MTM2016-78807-C2-1-P}

\maketitle
\begin{abstract}
As an application of the Bochner formula, we prove that if a $2$-dimensional Riemannian manifold admits a non-trivial smooth tangent vector field $X$ then its Gauss curvature is the divergence of a tangent vector field, constructed from $X$, defined on the open subset out the zeroes of $X$. Thanks to the Whitney embedding theorem and a standard approximation procedure, as a consequence, we give a new proof of the following well-known fact: if on an orientable, connected and compact $2$-dimensional smooth manifold there exists a continuous tangent vector field with no zeroes, then the manifold is diffeomorphic (or equivalently homeomorphic) to a torus.
\end{abstract}

\markboth{J.~M.~Almira, A. Romero}{Topology of orientable connected and compact surfaces}

\section{Introduction}
One of the most famous results in topology is the so-called  ``Hairy Ball Theorem'' which states that every continuous tangent vector field defined on an even-dimensional sphere vanishes somewhere, a fact that can be rephrased in several forms and has some interesting consequences. For example, the result can be restated by claiming that, if $n=2N$ is an even natural number, then there are no continuous unitary tangent vector fields on $\mathbb{S}^n$, and can be used to demonstrate that every continuous map $f:\mathbb{S}^{2N}\to \mathbb{S}^{2N}$ has either a fixed point or an antipodal point (a point $x$ such that $f(x)=-x$) \cite{Wi}. 

The best-known proof of the Hairy Ball Theorem is due to Brouwer \cite{Bro} and is based on the homology of the spheres.  Other proofs use degree theory \cite{E} and other algebraic topology arguments \cite{Cro, Cu}, but there are also analytic proofs \cite{Milnor1}, combinatorial proofs \cite{JT}, proofs based on the vector analysis of differential forms \cite{Bo1}, on differential topology \cite{Milnor, Poin},  etc. The proof we present in this paper is based on Riemannian geometry arguments, is restricted to the case $n=2$ and, as far as we know, it is the first one of this type that exists. 

The case $n=2$ has special relevance and was first proved by Poincaré \cite{Poin} using his index theorem, which characterizes Euler characteristic as the sum of the indices of a tangent vector field whose singularities are isolated. In particular,  if a nowhere vanishing tangent vector field exists on a connected compact surface  $M$, then $\chi(M)=0$, which represents an obstruction to the existence of such tangent vector fields on the sphere, since $\chi(\mathbb{S}^2)=2$. Note that Euler's characteristic can be used to reformulate the result as a characterization of the topological torus. Indeed, if $M$ is an orientable connected compact $2$-dimensional smooth manifold, then $\chi(M)=2-2g$, where $g$ is the genus of $M$. Thus, $M$ is a torus if and only if $\chi(M)=0$, since $g=1$ is the unique solution of $2-2g=0$. On the other hand, it is not difficult to construct nowhere vanishing  tangent vector fields on the torus (a fact that, surprisingly, has relevant applications in fusion theory \cite[Chapter 4]{paris}). Thus, it is natural to state Hairy Ball Theorem as follows:

\begin{theorem} \label{HBT}
 Assume $M$ is a $2$-dimensional orientable, connected and compact smooth manifold. Then the following are equivalent statements:
 \begin{itemize}
 \item[$(a)$] $M$ admits a continuous tangent  vector field with no zeroes.
 \item[$(b)$]  $M$  is diffeomorphic (or equivalently homeomorphic) to  a torus.
 \end{itemize}
\end{theorem}

In this note, we present a new proof of Theorem \ref{HBT}. The main ingredients of the proof are, after endowing $M$ with a Riemannian metric $g$, a combination of the Bochner technique, Gauss-Bonnet theorem, and a standard approximation argument. We only demonstrate the implication $(a)\Rightarrow (b)$, since the other implication is an easy exercise that consists of the explicit definition of a nowhere vanishing  tangent vector field on the torus.

\section{Proof of the main result}
\noindent One of the most important techniques in Riemannian geometry is the Bochner technique \cite[Chapter 9]{Pe}, \cite{W}. Roughly, it consists of relating the existence of certain non-trivial tangent vector fields on a Riemannian manifold with properties of the curvature of such a manifold. In the compact case, it leads to the vanishing of certain geometrically relevant tangent vector fields (Killing or conformal) under the assumption of positive or negative curvature everywhere \cite{W}. 
\vspace{3mm}

Our argument for the proof of $(a)\Rightarrow (b)$ in Theorem \ref{HBT} will split into three steps. In the first one, given a $2$-dimensional Riemannian manifold $(M,g)$, if it admits a smooth tangent vector field with no zeroes,  we use the well known Bochner formula (see, e.g., \cite{W,R} ) to construct a tangent vector field on $M$ whose divergence, with respect to $g$, is equal to the Gauss curvature of $g$. In the second step, we will start from a compact and connected $2$-dimensional smooth manifold $M$ which admits a nowhere zero continuous tangent vector field. We will endow $M$ with a Riemannian metric $g$. Then we will prove, by a classical approximation argument, that existence of the nowhere zero continuous tangent vector field leads to the existence of a smooth one also with no zeroes. Finally, in the last step, we will make use of the global Gauss-Bonnet theorem \cite{B,F} to get that $M$ must be a torus.

\vspace{3mm}

\noindent \textbf{Step 1. An application of the Bochner formula}

\noindent Consider a compact, orientable and connected $2$-dimensional manifold $M$. The classical Whitney embedding theorem \cite[Theorem 5.4.8]{Bo} guarantees that we can get a smooth embedding of $M$ into $\mathbb{R}^n$, for $n$ large enough. Then we can consider on $M$ the metric $g$ induced by the usual one of  $\mathbb{R}^n$ and this is precisely what we do. 
\vspace{3mm}

\noindent The Bochner formula states that
$$
X(\mathrm{div}(X))=-\mathrm{Ric}(X,X)+\mathrm{div}(\nabla_{X}X)-\mathrm{trace}(A_X^2),
$$
for all smooth tangent vector field $X$, where Ric is the Ricci tensor of $(M,g)$, div denotes the divergence on $(M,g)$ and $A_X$ is the linear operator field defined by $A_X(v)=-\nabla_vX$, being $\nabla$ the Levi-Civita connection of $g$ and $v\in T_xM$,  $x\in M$ \cite{W, R}. 

\vspace{1mm}

Now observe that $$\mathrm{div}(X)=-\mathrm{trace}(A_X),$$for any $X$ and any dimension of $M$. 

\vspace{1mm}

On the other hand, when  $\dim M=2$, we have $$\mathrm{Ric}(X,Y)=K\,g(X,Y),$$ for all $X,Y$,  where $K$ denotes the Gauss curvature of $g$. 

\vspace{1mm}

We want to apply Bochner formula to the unitary tangent vector field $\mathcal{T}:=T/g(T,T)^{\frac{1}{2}}$ in the case $\dim M=2$. 
Taking into account $g(\mathcal{T},\mathcal{T})=1$, we have that the matrix of $A_{\mathcal{T}}$ respect to an orthonormal basis $(\mathcal{T},E)$ is
\[
\begin{pmatrix} 
0 & 0 \\ a & b 
\end{pmatrix}.
\] 
\vspace{-2mm}

\noindent Therefore, we have $$\mathrm{trace}(A_{\mathcal{T}}^2)=(\mathrm{trace}(A_{\mathcal{T}}))^2=(\mathrm{div}(\mathcal{T}))^2.$$ 
Hence, we arrive to the following formula, 
\begin{equation} \label{bochner_formula}
\mathcal{T} (\mathrm{div} (\mathcal{T}))=-K+\mathrm{div}(\nabla_{\mathcal{T}}\mathcal{T})-(\mathrm{div}(\mathcal{T}))^2.
\end{equation}

On the other hand, using the well known formula 
$$\mathrm{div}(f\,X)=X(f)+f\mathrm{div}(X),$$
which holds true for any tangent vector field $X$ and any smooth function $f$ on $M$, we get
\begin{equation} \label{divergencia}
(\mathrm{div}(\mathcal{T}))^2= \mathrm{div}\big(\mathrm{div}(\mathcal{T})\,\mathcal{T}\big)-\mathcal{T}(\mathrm{div}(\mathcal{T})),
\end{equation}
and substituting (\ref{divergencia}) in \eqref{bochner_formula}, we get 
\begin{eqnarray*}
K &=&  \mathrm{div}(\nabla_{\mathcal{T}}\mathcal{T})-\mathrm{div}(\mathrm{div}(\mathcal{T})\,\mathcal{T})\\[2mm] 
 &=& \mathrm{div}\Big(\nabla_{\mathcal{T}}\mathcal{T}-\mathrm{div}(\mathcal{T})\,\mathcal{T}\Big).\\ 
\end{eqnarray*}

Thus, taking $Y=\nabla_{\mathcal{T}}\mathcal{T}-\mathrm{div}(\mathcal{T})\mathcal{T}$ we have the announced formula 

\begin{equation} \label{curvature}
K=\mathrm{div}(Y)\,.
\end{equation}

\vspace{0.5cm}

\noindent \textbf{Step 2. The approximation argument}

\vspace{2mm}
%

\vspace{1mm} 

Every continuous tangent vector field $X$ on $M$ may be identified with a continuous map $$X=(X_1,\cdots,X_n) \,:\, M\, \longrightarrow \,\mathbb{R}^n,$$
satisfying that $X(x)$ lies in the tangent space $T_x(M)$ (contemplated into $\mathbb{R}^n$ via the embedding), for all $x\in M$. Now, assume a continuous tangent vector field $X$ that vanishes nowhere and consider $$Z:=\frac{X}{\|X\|}=(Z_1,\cdots ,Z_n),$$ where $\|\cdot\|$ denotes the Euclidean norm of 
$\mathbb{R}^n$. Note that every function $Z_i:M\to\mathbb{R}$, $1\leq i \leq n$, is continuous. 

\vspace{1mm}

On the other hand, compactness of $M$ implies that $$M\,\subset\, \overline{\mathbf{B}}_0(r)=\{x\in\mathbb{R}^n:\|x\|\leq r\},$$ for some $r>0$. Thus, Tietze's extension theorem \cite{Man, O} implies that there exist functions $F_i\in\mathbf{C}(\overline{\mathbf{B}}_0(r),\mathbb{R})$ such that ${F_i}_{|M}=Z_i$ for $i=1,\cdots,n$. 

\vspace{1mm}

Apply now the Stone-Weierstrass approximation theorem \cite[Theorem 4.1]{Pi} (see also \cite{S2}) to the space of continuous functions $\mathbf{C}(\overline{\mathbf{B}}_0(r),\mathbb{R})$ and its subalgebra $\mathbb{R}[x_1,\cdots,x_n]$ to conclude that there are polynomials $P_i\in\mathbb{R}[x_1,\cdots,x_n]$, $1\leq 1 \leq  n$, such that 
$$\sup_{\|\mathbf{x}\|\leq r}|F_i(\mathbf{x})-P_i(\mathbf{x})| < \frac{1}{2\sqrt{n}}, \  i=1,\cdots,n,$$ 
which implies that 
\begin{eqnarray*}
\sup_{\mathbf{x}\in M}\left(\sum_{i=1}^n|Z_i(\mathbf{x})-P_i(\mathbf{x})|^2\right)^{\frac{1}{2}} &\leq&  \sup_{\|\mathbf{x}\|\leq r}\left(\sum_{i=1}^n|F_i(\mathbf{x})-P_i(\mathbf{x})|^2\right)^{\frac{1}{2}}\\ 
&\leq& \left(n (\max_{1\leq i\leq n}\sup_{\|\mathbf{x}\|\leq r}|F_i(\mathbf{x})-P_i(\mathbf{x})| )^2 \right)^{\frac{1}{2}}\\
&<& \left(n \frac{1}{4n} \right)^{\frac{1}{2}}=\sqrt{\frac{1}{4}} = \frac{1}{2}\,.
\end{eqnarray*}

\vspace{1mm}

In particular, $$P=(P_1,\cdots,P_n) : M \longrightarrow \mathbb{R}^n,$$ defines a smooth map. Let us now consider the tangential component $T$ of $P$. In other words, we consider, for each $\mathbf{x}\in M$, the decomposition of $T_\mathbf{x}\mathbb{R}^n$ as a direct sum $T_\mathbf{x}\mathbb{R}^n= T_\mathbf{x}M\oplus T_\mathbf{x}^{\perp}M$ and decompose $$P=T+U,$$ along the embedding, with $T(\mathbf{x})\in T_\mathbf{x}M$ and $U(\mathbf{x})\in T_\mathbf{x}^{\perp}M$, for each 
$\mathbf{x}\in M$. It follows that $T(\mathbf{x})\neq 0$ for all $\mathbf{x}\in M$ since $T(\mathbf{x}_0)=0$ for some 
$\mathbf{x_0}\in M$ implies that  
\[
\frac{1}{4}>\|P(\mathbf{x}_0)-Z(\mathbf{x}_0)\|^2=\|U(\mathbf{x}_0)-Z(\mathbf{x}_0)\|^2=\|U(\mathbf{x}_0)\|^2+1>1,
\]
which is impossible (note that in the formula above we have used that $U(\mathbf{x}_0)\perp Z(\mathbf{x}_0)$). Hence $T$ is a smooth tangent vector field with no zeroes on $M$.  

\vspace{3mm}

\noindent {\bf Remark.}
Obviously, the same conclusion can be obtained for $M$ of any dimension.

\vspace{3mm}

\noindent \textbf{Step 3. The Gauss-Bonnet argument}

\vspace{3mm}

\noindent If we apply the classical divergence theorem taking in mind (\ref{curvature}), we have that $$\int_MK\,dA=0,$$ (compare with \cite[Theorem p. 68]{F}). However, from the Gauss-Bonnet theorem we know $$\int_MK\,dA=2\pi \chi(M),$$ where $\chi(M)$ is the Euler number of $M$. Therefore, $\chi(M)=0$ and consequently $M$ is a torus \cite[Theorem 9.3.11]{H}.

\section{Acknowledgement} We want to express our acknowledgement to the referee, since his/her many suggestions have deeply improved the clarity of our arguments, making the paper more transparent and self-contained.

\bigskip

\begin{flushleft}

{\bf AMS Subject Classification: 53C20, 53C23, 57R25}\\[2ex]

{\bf Keywords:} Bochner technique, Orientable compact connected surfaces, Gauss curvature, Gauss-Bonnet theorem \\[2ex]

%
J. M.~ALMIRA,\\
Departamento de Ingenier\'{\i}a y Tecnolog\'{\i}a de Computadores,  Universidad de Murcia. \\
 30100 Murcia, SPAIN\\
e-mail: \texttt{jmalmira@um.es}\\[2ex]

A.~ROMERO,\\
Departamento de Geometr\'{\i}a y Topolog\'{\i}a, Universidad de Granada\\
    18071 Granada,  SPAIN.\\
e-mail: \texttt{aromero@ugr.es}\\[2ex]
%

\end{flushleft}


\end{document}